\documentclass[preprint]{iucr}              

     \journalcode{S}              

\usepackage{bm}
\usepackage{amsmath}
\usepackage{amsfonts}
\usepackage[english]{babel}
\usepackage[utf8]{inputenc}
\usepackage{algorithm}
\usepackage[noend]{algpseudocode}

\makeatletter
\newcommand{\mycaption}{%
\ifx \@captype \@undefined \@latex@error {\noexpand \caption outside float}\@ehd \expandafter \@gobble \else \refstepcounter \@captype \expandafter \@firstofone \fi {\@dblarg {\@caption \@captype }}%
}%
\makeatother

 \usepackage{color,soul}

\graphicspath{ {revised_figures/} }
\usepackage[colorinlistoftodos]{todonotes}

\begin{document}                  



\title{Improving reproducibility in synchrotron tomography using implementation-adapted filters} 
\shorttitle{Implementation-adapted filters}


\author[a,b]{Poulami Somanya}{Ganguly}
\author[a,c]{Dani\"{e}l M.}{Pelt}
\author[d,e]{Doga}{G\"{u}rsoy}
\author[d]{Francesco}{de Carlo}
\author[a,c]{K. Joost}{Batenburg}

\aff[a]{Computational Imaging, Centrum Wiskunde \& Informatica, Amsterdam, \country{The Netherlands}}
\aff[b]{Mathematical Institute, Leiden University, Leiden, \country{The Netherlands}}
\aff[c]{Leiden Institute of Advanced Computer Science, Leiden University, Leiden, \country{The Netherlands}}
\aff[d]{X-ray Science Division, Argonne National Laboratory, Argonne, IL, \country{USA}}
\aff[e]{Department of Electrical Engineering and Computer Science,
Northwestern University, Evanston, IL, \country{USA}}









\maketitle                        


\begin{synopsis} 
Dissimilar hardware and software conventions at various synchrotrons lead to quantitative differences in experimental results. This paper proposes a method to improve reproducibility of tomographic reconstructions by optimising the filtering step in commonly-used reconstruction algorithms.
\end{synopsis}

\begin{abstract} 


For reconstructing large tomographic datasets fast, filtered backprojection-type or Fourier-based algorithms are still the method of choice, as they have been for decades. These robust and computationally efficient algorithms have been integrated in a broad range of software packages. The continuous mathematical formulas used for image reconstruction in such algorithms are unambiguous. However, variations in discretisation and interpolation result in quantitative differences between reconstructed images, and corresponding segmentations, obtained from different software. This hinders reproducibility of experimental results, making it difficult to ensure that results and conclusions from experiments can be reproduced at different facilities or using different software.

In this paper, we propose a way to reduce such differences by optimising the filter used in analytical algorithms. These filters can be computed using a wrapper routine around a black-box implementation of a reconstruction algorithm, and lead to quantitatively similar reconstructions. 
We demonstrate use cases for our approach by computing implementation-adapted filters for several open-source implementations and applying it to simulated phantoms and real-world data acquired at the synchrotron. Our contribution to a reproducible reconstruction step forms a building block towards a fully reproducible synchrotron tomography data processing pipeline.
\end{abstract}


\section{Introduction}\label{intro}




In several scientific disciplines, such as materials science, biomedicine and engineering, a quantitative three-dimensional representation of a sample of interest is crucial for characterising and understanding the underlying system \cite{fusseis2014brief,luo2018cracking,midgley2009electron, rubin2014computed}. Such a representation can be obtained with the experimental technique of computerised tomography (CT). In this approach, a penetrating beam, such as X-rays, is used to obtain projection images of a sample at various angles. These projections are then combined by using a computational algorithm to give a 3D reconstruction \cite{buzug2011computed, kak2002principles}. 

Different tomographic setups are used in various practical settings. Our focus here is on tomography performed with a \emph{parallel-beam} X-ray source at synchrotrons. Synchrotrons provide a powerful source of X-rays for imaging, enabling a broad range of high-resolution and high-speed tomographic imaging techniques \cite{thompson1984computed, de2006x, stock2019microcomputed}.

\begin{figure}
    \centering
    \includegraphics[scale=0.4]{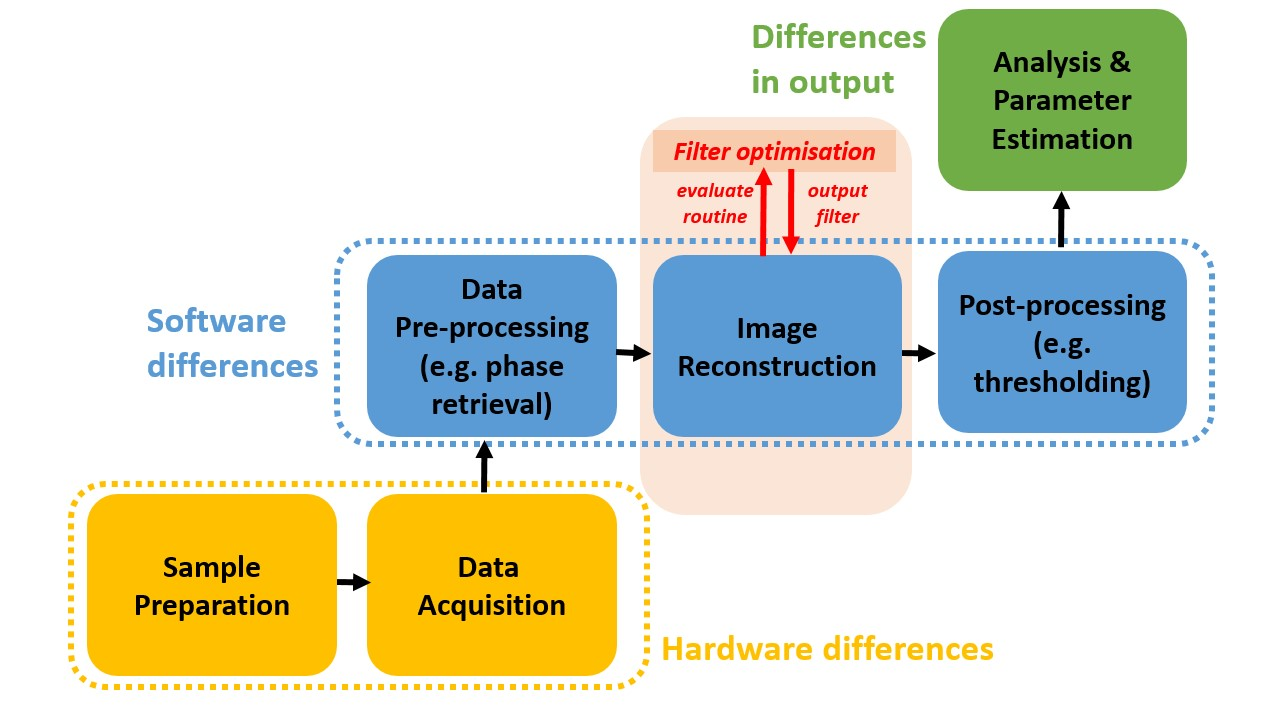}
    \caption{Schematic representation of a typical tomography pipeline at synchrotrons. Hardware differences play an important role during sample preparation and data acquisition. Software differences affect image pre-processing, reconstruction and post-processing. Together these lead to differences in the output of analysis and parameter estimation studies. In this paper we propose a filter optimisation method that works as a wrap-around routine on the reconstruction block. Our method only requires evaluations of the reconstruction routine and does not require any internal coding. The output of our method is a filter that can be used in the reconstruction block for more reproducible reconstructions.}
    \label{fig:pipeline}
\end{figure}

A typical tomography experiment at the synchrotron can be described by a pipeline consisting of several sequential steps (see Fig.~\ref{fig:pipeline}). First, a sample is prepared according to the experiment and imaging setup requirements. Then, the imaging system is aligned \cite{yang2017convolutional}, and a series of projection images of the sample are acquired \cite{hintermuller2010image}. These data are then processed for calibration, contrast improvement (e.g.~phase retrieval \cite{paganin2002simultaneous}) or removal of undesirable artefacts like rings or stripes \cite{massimi2018improved}. Following pre-processing, the data are fed into a reconstruction software package that makes use of one or more standard algorithms to compute a 3D reconstruction \cite{gursoy2014tomopy, pelt2016integration}. The reconstruction volumes can then be further post-processed and analysed \cite{salome1999synchrotron, buhrer2020unveiling} to obtain parameter estimates of the system being studied. In some cases, systematic imperfections in the data can also be corrected by post-processing reconstructions. For example, ring artefacts, which are commonly observed in synchrotron data, can be corrected before or after reconstruction \cite{gursoy2014tomopy}.

At various synchrotron facilities in the world, the pipeline described above is implemented using different instruments, protocols and methods specific for each facility \cite{kanitpanyacharoen2013comparative}. These differences are on the level of both hardware and software. Dissimilarities in the characteristics of the used X-ray source and detection system, including camera, visible light objective and scintillator screen, lead to differences in the acquired data. The differences in the data are then further compounded by variations in processing and reconstruction software, resulting in differences in voxel or pixel intensities, and eventually in variations in the output of post-processing and analysis routines.

For users, such differences pose several challenges. First, it is difficult to ensure that results and conclusions obtained from experiments at one facility are comparable and consistent with experiments from another facility. Second, other researchers seeking to reproduce the results of a previous work with their own software might not be able to do so, even if they have access to raw data. In \cite{kanitpanyacharoen2013comparative}, the authors report quantitative differences at various stages of the pipeline when scanning the same object at different synchrotrons. Reproducibility and the ability to verify experimental findings is crucial for ascertaining the reliability of scientific results. Therefore, in order to ensure reproducibility for the synchrotron pipeline, it is important to quantify and mitigate differences in the acquired, processed and reconstructed data.


Hardware and software vary across synchrotrons for a number of reasons.
Each synchrotron uses a pipeline that is optimised for its specific characteristics. In addition, legacy considerations play a role in the choice of components.
Because of the variations across synchrotrons, any successful strategy for creating reproducible results must take this diversity into account. Ideally, the choices for specific implementations of each block in the synchrotron pipeline in Fig~\ref{fig:pipeline} should not influence the final results of a tomography experiment. Following this strategy, each block can be optimised for reproducibility independently from the rest of the pipeline. 

In this paper, we focus on improving the reproducibility of the reconstruction block in the pipeline. In most synchrotrons, fast analytical methods such as filtered backprojection (FBP) \cite{kak2002principles} and Gridrec \cite{dowd1999developments} are the most commonly-used algorithms for reconstruction. This is primarily because
such algorithms are fast and work out-of-the-box without parameter tuning. These algorithms give accurate reconstructions when the projection data are well-sampled, such as in microCT beamlines where thousands of projections can be acquired in a relatively short time. 

Several open-source software packages for synchrotron tomography reconstruction are available, such as TomoPy, the ASTRA toolbox and scikit-image \cite{gursoy2014tomopy, palenstijn2013astra, van2014scikit}. Usually, an in-house implementation of FBP or Gridrec, or one of the open-source software packages is used for reconstruction. Each of these implementations contains a \emph{filtering} step
that is applied to the projection data as part of the reconstruction. Filtering influences characteristics, such as noise and smoothness, of the reconstructed volume. A sample-independent, pre-defined filter is generally used for reconstruction. Some filters used in this step have tunable parameters, but these are often tuned on-the-fly and are not recorded in metadata. 
    
Reconstructions in analytical algorithms are obtained by inversion of the Radon transform \cite{natterer2001mathematics}. Although this inversion is well-defined mathematically in a continuous setting, software implementations invariably have to work in a discretised space. In software implementations, the measurements as well as the reconstructed volume are \emph{discrete}. In a discretised space, inversion of the Radon transform often translates to a \emph{backprojection} step,
which makes use of a discretised \emph{projection kernel} to simulate the intersection between the scanned object and X rays \cite{pchansenbook2021}. The backprojection operation can also be performed directly using interpolations in Fourier space \cite{kak2002principles}. 

\begin{figure}
\includegraphics[scale=0.28]{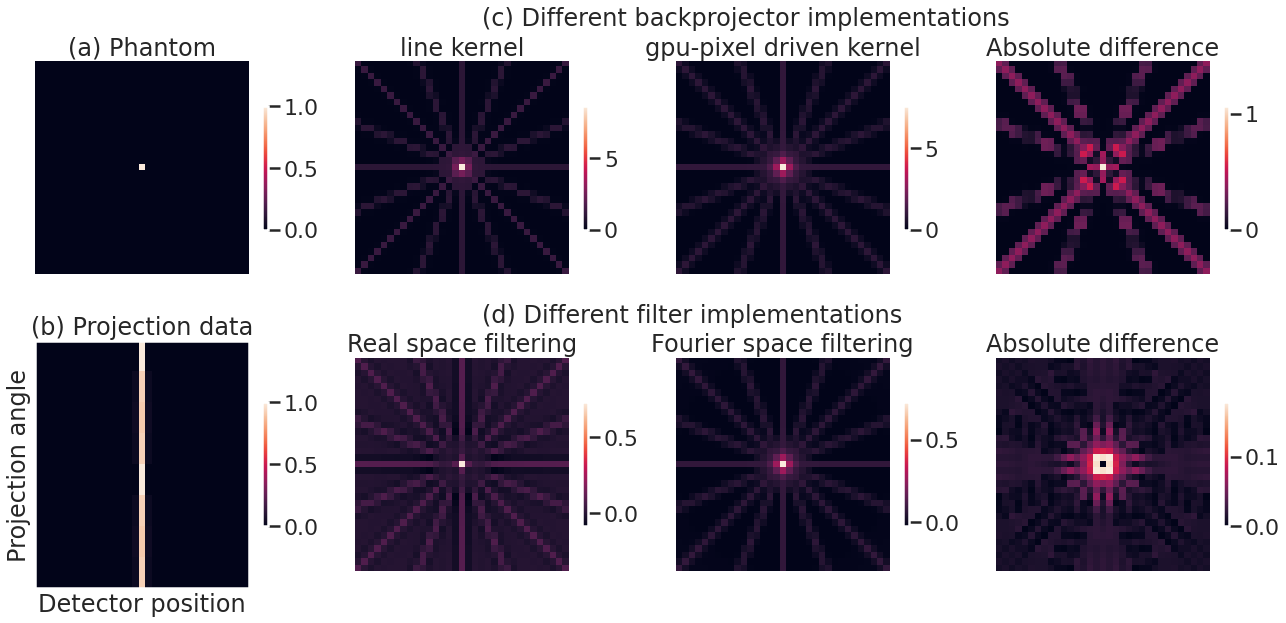}
\caption{Differences in reconstruction due to differences in backprojector and filter implementations. (a) a $33 \times 33$ phantom with one bright pixel, (b) sinogram of the phantom (computed using a strip kernel from the ASTRA toolbox), (c) differences in (unfiltered) backprojection when using different backprojectors: (\textit{left to right}) backprojection using a CPU line kernel from the ASTRA toolbox, backprojection using a GPU pixel-driven kernel from the ASTRA toolbox, absolute difference between the two backprojections. (d) differences in reconstruction when using different filtering routines in FBP with the \texttt{gpu-pixel} kernel as backprojector: (\textit{left to right}) reconstruction using filtering in real space with the Ram-Lak filter, reconstruction using the ramp filter in Fourier space, absolute difference between the two reconstructions.}
    \label{fig:diffs}
\end{figure}

Different choices of discretisation and interpolation, in projection kernels and filters, are possible. These choices lead to quantitative differences between the reconstructions obtained from different software implementations. A simple example of this effect is shown in Fig.~\ref{fig:diffs}, where we consider a phantom of pixel size $33 \times 33$ and data along $8$ projection angles uniformly sampled in $[0,\pi)$. We compare reconstructions of the same data using two different projection kernels and two different filtering methods. In both instances, the image to be reconstructed contains a single bright pixel at the centre of the field-of-view. The \emph{sinogram} of such an image (i.e. the combined projection data for the full range of angles) was computed using a CPU strip kernel projector from the ASTRA toolbox \cite{palenstijn2013astra}. Backprojections of this projection data using two other projectors - a CPU line projection kernel and a pixel-driven kernel implemented on a graphics processing unit (GPU) - show significant, radially-symmetric differences. These differences are dependent on the number of projection angles used, and are highly structured, unlike differences due to random noise. We also observe structured differences between reconstructions when the same projection kernel (\texttt{gpu-pixel}) is used after different filtering operations in real and Fourier space.
This example highlights the impact of discretisation and interpolation choices on the final reconstruction obtained from identical raw data. 

Our main contribution in this paper is a heuristic approach that can improve reproducibility in reconstructions.
Our method consists of optimising the filter used in different software implementations of reconstruction methods. We call such optimised filters \emph{implementation-adapted filters}. The computation of our filters does not require knowledge of the underlying software implementation of the reconstruction algorithm. Instead, a wrapper routine around any black-box implementation can be used for filter computation. Once computed, these filters can be applied with the reconstruction software like any other standard filter. 

Our paper is organised as follows. In Section \ref{background}, we formulate the reconstruction problem mathematically and discuss the effect of different software implementations. In Section \ref{filters}, we describe our algorithm for computing implementation-adapted filters. Numerical experiments described in Sections \ref{data_and_metrics} and \ref{results} demonstrate use cases for our filters on simulated and real data. Finally, we discuss extensions to the current work in Section \ref{discussion} and conclude our paper in Section \ref{conclusion}. Our open-source Python code for computing implementation-adapted filters is available on GitHub\footnote{{https://github.com/poulamisganguly/impl-adapted-filters}}.

\section{Background}\label{background}

\subsection{Continuous reconstruction}

Consider an object described by a two-dimensional attenuation function $f: \mathbb{R}^2 \rightarrow \mathbb{R}$. Mathematically, the tomographic projections of the object can be modelled by the Radon transform, $\mathcal{R}(f)$. The Radon transform is the line integral of $f$ along parametrised lines $l_{\theta,t}= \{(x,y) \in \mathbb{R}^2 \,|\, x\cos\theta + y\sin\theta = t \}$, where $\theta$ is the projection angle and $t$ is the distance along the detector. Projection data $p_\theta(t)$ along an angle $\theta$ are thus given by
\begin{equation}\label{eq:cont_radon}
    p_\theta(t) = \mathcal{R}(f) = \iint_{\mathbb{R}^2} f(x, y) \delta(x\cos\theta + y\sin\theta - t) dx dy.
\end{equation}

The goal of tomographic reconstruction is to obtain the function $f(x,y)$ given the projections $p_\theta(t)$ for various angles $\theta \in \Theta$. One way to achieve this is by direct inversion of the Radon transform. Given a complete angular sampling in $[0,\pi)$, the Radon transform can be inverted giving the following relation \cite{kak2002principles}
\begin{equation}\label{eq:radon_inversion}
    f(x,y) = \int_0^\pi \Bigg( \int_{-\infty}^{\infty} \tilde{P}_\theta(\omega) |\omega| e^{2 \pi i \omega (x \cos \theta + y \sin \theta)} d\omega \Bigg) d\theta,
\end{equation}
where $\tilde{P}_\theta(\omega)$ denotes the Fourier transform of the projection data $p_\theta(t)$ and multiplication by the absolute value of the frequency $|\omega|$ denotes filtering with the so-called ramp filter.

For noiseless and complete data, the Radon inversion formula \eqref{eq:radon_inversion} provides a perfect analytical reconstruction of the function $f(x,y)$ from its projections. However, in practice, tomographic projections are obtained on a \emph{discretised} detector, consisting of individual pixels, and for a finite set of projection angles. Additionally, the reconstruction volume must be discretised in order to represent it on a computer. Therefore, in practical applications, a discretised version of \eqref{eq:radon_inversion} is used to obtain reconstructions.

\subsection{Discrete reconstruction}
Discretisation of the reconstruction problem yields the following equation for the discrete reconstruction $r(x_d, y_d)$:
\begin{equation}\label{eq:fbp}
    r(x_d, y_d) = \sum_{\theta_d \in \Theta}\sum_{t_d \in T} h(t_d) P_{\theta_d}(x_d \cos \theta_d + y_d \sin \theta_d - t_d),
\end{equation}
where $(x_d,y_d)$, $\theta_d$ and $t_d$ denote discretised reconstruction pixels, angles and detector positions, respectively, and $h(t_d)$ is a discrete real-space filter. This inversion formula is known as the filtered backprojection (FBP) algorithm.

The FBP equation \eqref{eq:fbp} can be written algebraically as the composition of two matrix operations: filtering and backprojection. Filtering denotes convolution in real space (or, correspondingly,  multiplication in Fourier space) with a discrete filter. Backprojection consists of a series of interpolation and numerical integration steps to sum contributions from different projection angles. These discretised operations can be implemented in a number of different ways and different software implementations often make use of different choices for discretisation and interpolation. Consequently, the reconstruction obtained from a particular implementation is dependent on these choices. The reconstruction $\bm{r}_I$ from an implementation $I$ can thus be written as
\begin{equation}\label{eq:analytical_reco}
    \bm{r}_I(\bm{h}, \bm{p}) = \bm{W}_I^T \bm{M}_I(\bm{h}, \bm{p}),
\end{equation}
where $\bm{W}_I^T$ is the backprojector and $\bm{M}_I(\cdot,\cdot)$ is the (linear) filtering operation associated with implementation $I$. We denote the discrete filter by $\bm{h}$.

In the following subsection, we discuss some common choices for projection and filtering operators in software implementations of analytical algorithms.

\subsection{Differences in projectors and filtering}\label{bp_operators}


In order to discretise the Radon transform, we must choose a suitable discretisation of the reconstruction volume, a discretisation of the incoming ray and an appropriate numerical integration scheme. All these choices contribute to differences in different backprojectors $\bm{W}_I^T$ in \eqref{eq:analytical_reco}. 

Voxels (or pixels in 2D) in the reconstruction volume can be considered either to have a finite size or to be spikes of infinitesimal size. Similarly, a ray can be discretised to have finite width (i.e.~a strip) or have zero width (i.e.~a line). The numerical integration scheme chosen might be piecewise constant, piecewise linear or continuous. All of these different choices have given rise to different software implementations of backprojectors \cite{pchansenbook2021}. 
There exist different categorisations of backprojectors in the literature; for example, the \texttt{linear} kernel in the ASTRA toolbox is referred to as the slice-interpolated scheme in \cite{xu2006comparative} and the \texttt{strip} kernel is referred to as the box-beam integrated scheme in the same work. In this paper, we designate different backprojectors with the terms used in the software package where they have been implemented. 

In addition to the choices mentioned above, backprojectors have also been optimised for the processing units on which they are used. For this reason, backprojectors that are optimised to be implemented on graphics processing units (GPUs) might be different from those that are implemented on a CPU due to speed considerations. In particular, GPUs provide hardware interpolation that is extremely fast, but can also be of limited accuracy compared to standard floating point operations.

So far, we have discussed real space backprojectors. Fourier-domain algorithms such as Gridrec \cite{dowd1999developments} use backprojectors that operate in the Fourier domain. These operators are generally faster than real-space operators, and are therefore particularly suited for accelerating iterative algorithms \cite{arcadu2016fast}. Unlike real space backprojectors, Fourier-space backprojectors perform interpolation in the Fourier domain. As this might lead to non-local errors in the reconstruction, an additional filtering step is performed to improve the accuracy of the interpolation.

Apart from differences in backprojectors, different implementations also vary in the way they perform the filtering operation in analytical algorithms. Filtering can be performed as a convolution in real space or as a multiplication in Fourier space. Real space filtering implementations can differ from each other in computational conventions, for example by the type of padding used \cite{marone2012regridding} to extend the signal at the boundary of the detector. Moreover, the zero-frequency filter component is treated in different ways between implementations. For example, the Gridrec implementation in TomoPy sets the zero-frequency component of the filter to zero. 



\section{Implementation-adapted filters}\label{filters}

We now present the main contribution of our paper. In order to mitigate the differences between implementations discussed in the previous section, we propose to specifically tune the filter $\bm{h}$ for each implemented analytical algorithm. In the following, we describe an optimisation scheme for the filter, which helps us to reduce the differences between reconstructions from various implementations.

We optimise the filter by minimising the $\ell^2$ difference with respect to the projection data $\bm{p}$. This can be stated as the following optimisation problem over filters $\bm{h}$:
\begin{equation}\label{eq:filter_opt}
    \bm{h}_I^\ast = \arg \min_{\bm{h}} \|\bm{p} - \bm{W}  \bm{r}_I(\bm{h}, \bm{p})\|^2_2,
\end{equation}
where $\bm{r}_I$ is the reconstruction from implementation $I$. Note that the forward projector $\bm{W}$ used above is chosen as a fixed operator in our method (the same for each implementation for which the filter is optimised) and does not have to be the transpose of the implementation-specific backprojection operator $\bm{W}_I^T$. In order to improve stability and take additional prior knowledge of the scanned object into account, a regularisation term can be added to the objective in \eqref{eq:filter_opt}.


The solution to the optimisation problem above is a implementation-adapted filter $\bm{h}_I^\ast$. Once the filter has been computed, it can be used in \eqref{eq:analytical_reco} to give an optimised reconstruction:
\begin{equation*}
    \bm{r}_I^\ast = \bm{W}_I^T \bm{M}_I(\bm{h}_I^\ast, \bm{p}).
\end{equation*}
Out of all reconstructions that an implemented algorithm can produce for a given dataset $\bm{p}$ by varying the filter, this reconstruction, $\bm{r}_I^\ast$, is the one that results in the smallest residual error. Such filters are known as minimum-residual filters and have previously been proposed to improve reconstructions of real-space analytical algorithms in low-dose settings \cite{pelt2014improving, lagerwerf2020automated}. 

Our implementation-adapted filters are thus minimum-residual filters that have been optimised to each implementation $I$. The main difference between the previous works \cite{pelt2014improving, lagerwerf2020automated} and our present study is that we use a fixed forward operator in our optimisation problem, which is not necessarily the transpose of the backprojection operator. More importantly, our goal in this paper is not the improvement of reconstruction accuracy, but the reduction of differences in reconstruction between various software implementations.

We hypothesise that such minimum-residual reconstructions obtained using different implementations are closer (quantitatively more similar) to each other than reconstructions obtained using standard filters. 
As an example for motivating this choice, let's take an implementation of an analytical algorithm from both TomoPy and the ASTRA toolbox. Given a certain dataset, changing the reconstruction filter results in different reconstructed images, each with a different residual error. Even though the implementations used by TomoPy and ASTRA are fixed, the freedom in choosing a filter gives us an opportunity to reduce the difference between reconstructions from both implementations.
Tuning the filter is a way to \emph{optimise} the reconstruction according to user-selected quality criteria. Choosing the \emph{minimum-residual} reconstruction for each implementation results in reconstructions that are the \emph{closest possible} to each other in terms of data misfit. Closeness in data misfit, under convexity assumptions, indicates closeness in pixel intensity values of reconstruction images. Hence, the minimum-residual reconstructions for the two implementations are closer to each other than reconstructions with standard filters offered by the implementations. 



To compute the optimised filter \eqref{eq:filter_opt}, 
we use the fact that the reconstruction $\bm{r}_I(\bm{h}, \bm{p})$ of data $\bm{p}$ obtained from an implementation of FBP or Gridrec is \emph{linear} in the filter $\bm{h}$. This means that we can write the reconstruction as
\begin{equation*}
    \bm{r}_I(\bm{h}, \bm{p}) = \bm{R}_I(\bm{p}) \bm{h}, 
\end{equation*}
where $\bm{R}_I(\bm{p})$ is the reconstruction matrix of implementation $I$ given projection data $\bm{p}$. Thus, the optimisation problem \eqref{eq:filter_opt} becomes
\begin{equation}\label{eq:filter_opt2}
    \bm{h}_I^\ast = \arg \min_{\bm{h}} \|\bm{p} - \bm{W}   \bm{R}_I(\bm{p}) \bm{h}\|^2_2 =: \arg \min_{\bm{h}} \|\bm{p} - \bm{F}_I(\bm{p}) \bm{h}\|^2_2
\end{equation}
The matrix $\bm{F}_I(\bm{p})$ has dimensions $N_p \times N_f$, where $N_p$ is the size of projection data and $N_f$ is the number of filter components. For a filter that is independent of projection angle, the number of filter components, $N_f$, is equal to the number of discrete detector pixels, $N_d$. The projection size $N_p := N_d N_\theta$, where $N_\theta$ is the number of projection angles. $\bm{F}_I(\bm{p})$ can be constructed explicitly by assuming a basis for filter components. A canonical basis can be formed using $N_d$ unit vectors $\{\bm{e}_i, i=1, 2, \ldots, N_d\}$, such that
\begin{equation*}
\bm{e}_1 = \begin{pmatrix}1\\0\\.\\.\\.\\0 \end{pmatrix},\quad
\bm{e}_2 = \begin{pmatrix}0\\1\\.\\.\\.\\0 \end{pmatrix}, \quad \ldots \quad
\bm{e}_{N_d} = \begin{pmatrix}0\\0\\.\\.\\.\\1 \end{pmatrix}.
\end{equation*}
Using these basis filters, each column of $\bm{F}_I(\bm{p})$ can be computed by reconstructing $\bm{p}$ using the implementation $I$, followed by forward projection with $\bm{W}$: 
\begin{eqnarray*}
    \bm{f}_j = \bm{W}\bm{r}_I(\bm{e}_j,\bm{p}), \quad j \in \{1,2,\ldots,N_d \}\\
    \bm{F}_I(\bm{p}) = \begin{pmatrix}
    \bm{f}_1 & \bm{f}_2 & \bm{f}_3 & \ldots & \bm{f}_{N_{d}}
    \end{pmatrix}
\end{eqnarray*}
We can then substitute for $\bm{F}_I(\bm{p})$ in \eqref{eq:filter_opt2} and solve for the optimised filter $\bm{h}^\ast_I$. Note that our method only requires \emph{evaluations} of the implementation $I$ by using it as a black-box routine to compute the reconstructions $\bm{r}_I(\bm{e}_j,\bm{p})$ above. In other words, no knowledge of the implementation $I$ or any internal coding is required. 


If we expand the filter in a basis of unit vectors, $\mathcal{O}(N_p)$ reconstructions using the implementation $I$ and $\mathcal{O}(N_p)$ forward projections with $\bm{W}$ must be performed for filter optimisation. In contrast, the complexity of a standard FBP reconstruction is of the order of a single backprojection. Choosing a smaller set of suitable basis functions would result in a reduction in the number of operations for filter optimisation and, consequently, faster filter computations. One way to do this is by exponential binning \cite{pelt2014improving}. 

The idea of exponential binning is to assume that the real-space filter is a piecewise constant function with $N_b$ bins, where $N_{b} < N_{d}$. The bin width $w_i, \text{ for } i=1,2,\ldots,N_b$, is assumed to increase in an exponential fashion away from the centre of the detector, such that:
\begin{equation}
    w_i = \begin{cases}
    1, & |i|<N_l\\
    2^{|i| - N_l}, & |i| \geq N_l
    \end{cases},
\end{equation}
where $N_l$ is the number of large bins with width $1$.
Exponential binning is inspired by the observation that standard filters used in tomographic reconstruction, such as the Ram-Lak filter, are peaked at the centre of the detector and decay to zero relatively quickly towards the edges. Binning results in a reduction of free filter components from $N_d$ to $N_b$. Moreover, despite the reduction in components, it does not typically result in a significant change in reconstruction quality \cite{pelt2014improving}. 

The pseudocode for our filter computation method is shown in Algorithm \ref{alg:alg_filter}. Here we give further details of the routines used in the algorithm. The \texttt{filter} routine performs filtering in the Fourier domain, which is equivalent to multiplication by the filter followed by an inverse Fourier transform. The $\texttt{reconstruct}_I$ routine calls the function for reconstruction in implementation $I$ with the internal filtering disabled. Finally, the \texttt{lstsq} routine calls a standard linear least squares solver in NumPy \cite{harris2020array} to compute filter coefficients. 

\begin{algorithm}
\mycaption{Implementation-adapted filter computation}\label{alg:alg_filter}
\begin{algorithmic}[1]
\Procedure{Compute filter}{$\bm{p}$, $I$, $\bm{W}$}:
\State Create filter basis: $\mathcal{B} := \{b_1,b_2,\ldots,b_{N_b}\}$
\State Compute columns of $\bm{F}_I(\bm{p})$:
\For{$\bm{b}_j \in \mathcal{B}$}
\State Filter data with basis filter: $\bm{q} \leftarrow \texttt{filter}(\bm{p},\bm{b}_j)$
\State Reconstruct filtered projection with $I$: $\bm{r} \leftarrow \texttt{reconstruct}_I(\bm{q})$
\State Forward project reconstruction $\bm{f}_j \leftarrow \texttt{flatten}(\bm{W}\bm{r})$
\EndFor
\State Linear least squares fitting of filter coefficients: $\bm{c} \leftarrow \texttt{lstsq}(\bm{F}_I(\bm{p}), \bm{p})$
\State Return filter: $\bm{h}^\ast \leftarrow \sum_{j=1}^{N_b} c_j \bm{b}_j$
\EndProcedure
\end{algorithmic}
\end{algorithm}

Once a filter $\bm{h}^\ast$ is computed, we can store it in memory, either as a filter in Fourier space or as a filter in real space after computing the Fourier transform of $\bm{h}^\ast$. Using the filter with a black-box software package involves calling the \texttt{filter} routine with the data and the computed filter as arguments, followed by one call of the $\texttt{reconstruct}_I$ routine in a chosen algorithm (with its internal filtering disabled). Thus, the complexity of a reconstruction using a computed implementation-adapted filter is the same as that of a reconstruction run using a standard filter.

In the following sections, we describe numerical experiments and the results of filter optimisation on reconstructions.


\section{Data and metrics}\label{data_and_metrics}

We performed a range of numerical experiments on real and simulated data to quantitatively assess (i) the effect of our proposed optimized filters on the variations between reconstructions from different implementations; (ii) the behaviour and dependence of our proposed filters on acquisition characteristics such as noise and sparse angular sampling; and (iii) the effect of our proposed filters on post-processing steps following the reconstruction block in Fig~\ref{fig:pipeline}. 
In this section, we describe the software implementations used, data generation steps and the metric used to quantify intra-set variability of reconstructions.

\subsection{Software implementations of analytical algorithms}\label{implementations}
We optimised filters to commonly used software implementations of FBP and Gridrec. For FBP, we considered different projector implementations in the ASTRA toolbox \cite{palenstijn2013astra} as well as the \texttt{iradon} backprojection function in scikit-image \cite{van2014scikit}. These implementations use different choices of volume and ray discretisation as well as numerical integration schemes. From the ASTRA toolbox, we considered projectors implemented on the CPU (\texttt{strip}, \texttt{line} and \texttt{linear}) as well as a pixel-driven kernel on the GPU (\texttt{gpu-pixel}, called \texttt{cuda} in the ASTRA toolbox). For Fourier-space methods, we considered the Gridrec implementation in TomoPy. We used the ASTRA \texttt{strip} kernel as the forward projector $\boldsymbol{W}$ in \eqref{eq:filter_opt} during filter computations.

\subsection{Projection data}
We performed experiments with both simulated and real data. Both data consisted of projections acquired in a parallel-beam geometry along a complete angular range in $[0, \pi)$.

\subsubsection{Simulated foam phantom data}
Simulated data of foam-like phantoms were generated using the foam\_ct\_phantom package in Python. This package generates 3D volumes of foam-like phantoms by removing, at random, a pre-specified number of non-overlapping spheres from a cylinder of a given material \cite{pelt2018improving}. The simulated phantoms are representative of real foam samples used in tomographic experiments and are challenging to reconstruct due to the presence of features at different length scales. At the same time, the phantoms are amenable to experimentation as data in different acquisition settings can be easily generated. Slices of one such phantom, which we used for the experiments in this paper, are shown in Fig.~\ref{fig:foam} and Fig.~\ref{fig:filter_variability}. 

Ray tracing through the volume is used to generate projection data from a 3D foam phantom. To simulate real-world experimental setups, where detector pixels have a finite area, ray supersampling can be used. This amounts to averaging the contribution of $n$ neighbouring rays within a single pixel, where $n$ is called the supersampling factor.

For our experiments, we generated a 3D foam with $1000$ non-overlapping spheres with varying radii. 
A parallel beam projection geometry, in line with synchrotron setups, was used to generate projection data. We used ray supersampling with a supersampling factor of 4, and each 2D projection was discretised on a pixel grid of size $256 \times 256$. We varied the number of projection angles, $N_\theta$, in our experiments in order to determine the effect of sparse sampling ranges on our filters.

Poisson noise was added to noiseless data by using the \texttt{astra.add\_noise\_to\_sino} function in the ASTRA toolbox \cite{palenstijn2013astra}. This function requires the user to specify a value for the photon flux $I_0$. In an image corrupted with Poisson noise, each pixel intensity value $k$ is drawn from a Poisson distribution
\begin{equation*}
   f_{\text{Pois}}(k;\lambda) = \dfrac{\lambda^k e^{- \lambda}}{k!},
\end{equation*}
with $\lambda \propto I_0$. High photon counts (and high values of $\lambda$) correspond to low noise settings. All noise realisations in our experiments were generated with a pre-specified random seed.

\subsubsection{Real data of shale}
In order to validate the applicability of our method to real data, we performed numerical experiments using microCT data of the Round-Robin shale sample N1 from the tomographic data repository Tomobank \cite{de2018tomobank}. We used data acquired at the Advanced Photon Source (APS) for our experiments. The Round-Robin datasets were acquired for characterising the porosity and microstructures of shale, and the same sample has been imaged at different synchrotrons (using the same experimental settings) for comparison of results \cite{kanitpanyacharoen2013comparative}. The dataset we used was acquired with a $10$x objective lens and had an effective pixel size of approximately $0.7 \mu$m. Each projection in the dataset had pixel dimensions $2048 \times 2048$, and data were acquired over $1500$ projection angles.
In order to simulate sparse angular range settings, we removed projections at intervals of $m = 2, 3, 4, 5 \text{ and } 10$ from the complete data.


\subsection{Quantitative metrics}
Reconstructions of a 3D volume from parallel beam data can be done slice-wise, because data in different slices (along the rotation axis) are independent of each other in a parallel beam geometry. Therefore, all our quantitative metrics were computed on individual slices. Reconstructed slices of the simulated foam phantom were discretised on a pixel grid of size $256 \times 256$. Reconstruction slices of the Round-Robin dataset were discretised on a pixel grid of size $2048 \times 2048$. All CPU reconstructions were performed on an Intel(R) Core(TM) i7-8700K CPU with 12 cores. GPU reconstructions were performed on a single Nvidia GeForce GTX 1070 Ti GPU with CUDA version 10.0.


We were interested in comparing the similarity between reconstructions in a \emph{set} of images, without having a reference reconstruction. We quantified the intra-set variability between reconstruction slices obtained from different implementations using the pixelwise standard deviation between these. For a set of reconstruction slices $\{ \bm{r}_I, I \in \mathcal{I}\}$ obtained using different implementations $I$, the standard deviation of a pixel $j$ is given by:
\begin{equation}\label{eq:pixelwise_std}
    {\sigma}_{j} = \sqrt{\dfrac{1}{N_I}\sum_{I \in \mathcal{I}}\Big((r_{I})_{j} - \bar{r}_{j}\Big)^2}; \qquad \bar{r}_{j} = \dfrac{1}{N_I} \sum_{I \in \mathcal{I}} (r_{I})_{j},
\end{equation}
where $(r_{I})_{j}$ is the intensity value of pixel $j$ in reconstruction $\bm{r}_I$ and $N_I$ is the total number of implementations. 

In our experiments, we reconstructed the same data using our set of implementations $\{I \in \mathcal{I}\}$, by using the Ram-Lak filter and the Shepp-Logan filter as defined in different packages, and then by using filters $\{\bm{h}^\ast_I, I \in \mathcal{I}\}$ \eqref{eq:filter_opt} that were optimised to those implementations. As a result, we achieved three sets of reconstructions: one set using the Ram-Lak filter, a second set using the Shepp-Logan filter and a third set using the implementation-adapted filters. We computed the pixelwise standard deviation \eqref{eq:pixelwise_std} over slices for all sets.

The mean standard deviation of a slice $S$ (with dimensions $N \times N$) is defined as the mean of pixelwise standard deviations in that slice:
\begin{equation}\label{eq:mean_std_dev}
    \bar{\sigma}^S = \dfrac{1}{N^2} \sum_{j \in J^S} {\sigma}_{j},
\end{equation}
where $J^S$ is the list of pixels in slice $S$.

In addition to the mean, the histogram of standard deviations \eqref{eq:pixelwise_std} provides important information about the distribution of standard deviation values in a slice. The \emph{mode} of this histogram is the value of standard deviation that occurs most, and the tail of the histogram indicates the number of large standard deviations observed. For reconstructions that are more similar to each other, we would expect the histogram to be peaked at a value close to $0$ and have a small tail.

In order to quantify the difference between a reconstruction slice and the ground truth (in experiments where a ground truth was available), we used the root mean squared error (RMSE) given by
\begin{equation}
    \text{RMSE}(\bm{r}_I) = \sqrt{\frac{1}{N^2}\sum(\bm{r}_I - \bm{r}_{gt})^2},
\end{equation}
where $\bm{r}_{gt}$ is the ground truth reconstruction. For a set of reconstructions we used the squared bias defined below to quantify the difference with respect to the ground truth:
\begin{equation}\label{eq:pixelwise_bias}
 \Big(\text{bias}(\{\bm{r}_I, I \in \mathcal{I}\}) \Big)^2 = \Big(\bar{\bm{r}} - \bm{r}_{gt}\Big)^2,
\end{equation}
where $\bar{\bm{r}} := \sum_{I \in \mathcal{I}} \frac{1}{N_I} \bm{r}_I$ is the mean over the set of reconstructions. The squared bias, similar to the standard deviation in \eqref{eq:pixelwise_std} is a pixelwise measure. The mean squared bias over a slice $S$ is obtained by taking the mean of \eqref{eq:pixelwise_bias} over all pixels in the slice.

In our experiments, we also quantify the effect of filter optimisation on later post-processing steps after reconstruction. To do this, we threshold a set of reconstructions using Otsu's method \cite{otsu1979threshold}, which picks a single threshold to maximise the variance in intensity between binary classes. To quantify the accuracy of the resulting segmentations and to compare the similarity in a set we used two standard metrics for segmentation analysis: the $F_1$ score and the Jaccard index. The $F_1$ score takes into account false positives (fp), true positives (tp) and false negatives (fn) in binary segmentation and is given by:
\begin{equation}\label{eq:f1_score}
    F_1 = \dfrac{\text{tp}}{\text{tp} + \frac{1}{2} (\text{fp} + \text{fn})}.
\end{equation}
The Jaccard index is the ratio between the intersection and union of two sets A and B. In our case, one set is the segmented binary image and the other set is the binary ground truth image:
\begin{equation}\label{eq:jaccard}
    J(A,B) = \dfrac{|A \cap B|}{|A \cup B|}.
\end{equation}

\section{Numerical experiments and results}\label{results}
In this section, we give details of our numerical experiments and discuss their results.
\subsection{Foam phantom data}
\subsubsection{Reduction in differences between reconstructions}
\begin{figure}
    \includegraphics[scale=0.2]{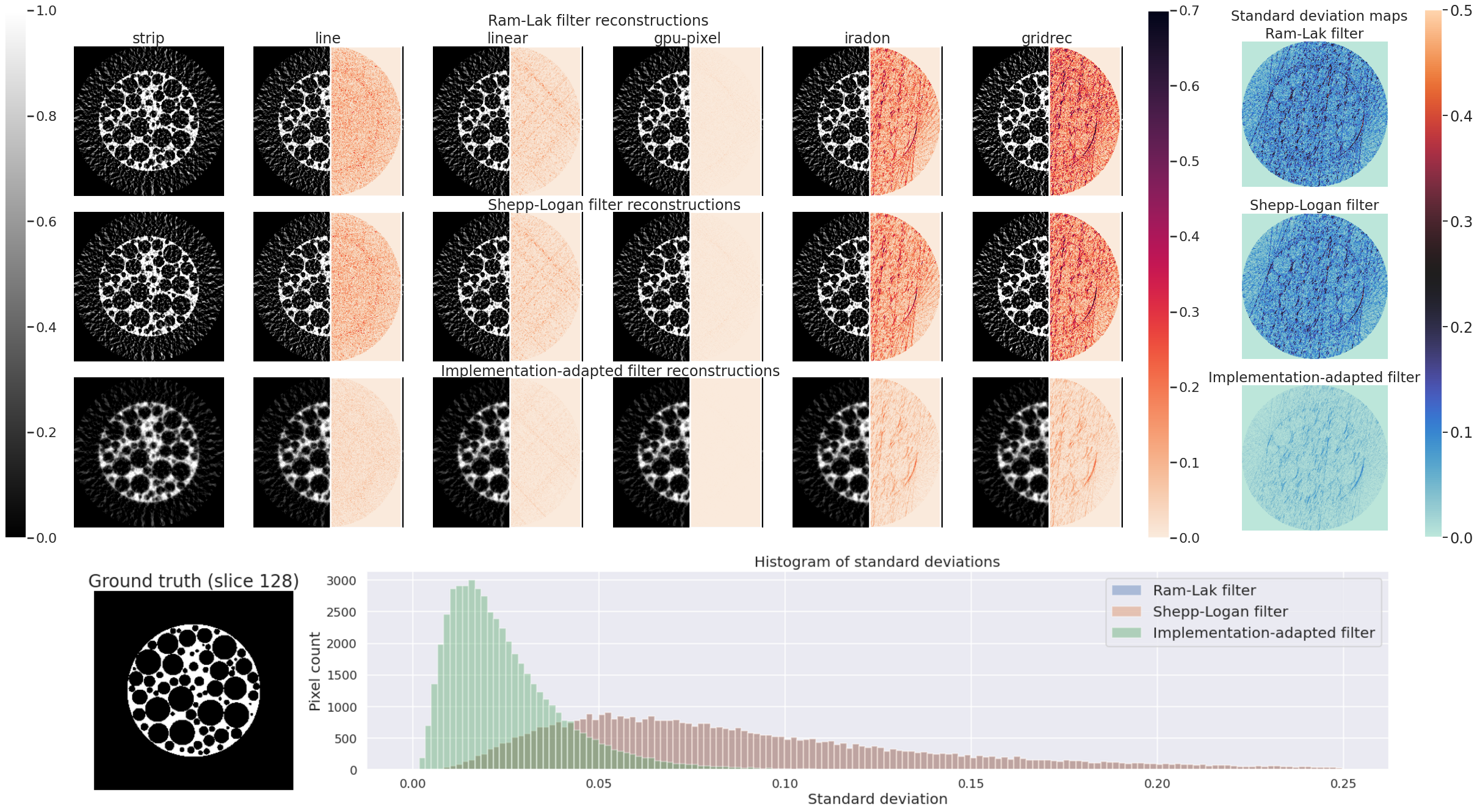}
    \caption{Reduction in intra-set variability between reconstructions of simulated foam data ($N_\theta = 32$, no noise) by using implementation-adapted filters. (\textit{top three rows}) Reconstructions of the central slice (slice no.~128) of a foam phantom. To highlight intra-set discrepancies we show the absolute difference with respect to the corresponding \texttt{strip} kernel reconstructions in the right half of each image. The rightmost column shows pixelwise standard deviation $\sigma$ in each set.
    (\textit{bottom row, left})
    Ground truth foam phantom slice. (\textit{right}) Histograms of standard deviations $\sigma$ for all three sets. The Ram-Lak filter and Shepp-Logan filter histograms overlap.}
    \label{fig:foam}
\end{figure}

Fig.~\ref{fig:foam} shows the central (ground truth) slice of the foam phantom. Data along $N_\theta = 32$ angles were reconstructed using all implementations using the Ram-Lak filter, the Shepp-Logan filter and our implementation-adapted filters. Reconstructions using the various filters are shown in Fig.~\ref{fig:foam}. In order to highlight intra-set variability, we include heatmaps showing the absolute difference with respect to one (\texttt{strip}) reconstruction. Upon visual inspection, we see that discrepancies between reconstructions are smaller in the set obtained using implementation-adapted filters. An interesting point to note is that the Gridrec and \texttt{iradon} reconstructions show the largest differences from the ASTRA \texttt{strip} kernel reconstruction in both sets. This suggests that differences between different software packages are greater than differences between different projectors in the same software package. 

To further investigate intra-set variability, we use pixelwise standard deviation maps for all sets of reconstructions. We observe that higher values of standard deviation are observed when using the Ram-Lak and Shepp-Logan filters. This indicates that quantitative differences between these reconstructions were more pronounced. In contrast, reconstructions using our implementation-adapted filters were more similar, resulting in low pixelwise standard deviations. Furthermore, the mode of the histogram of standard deviations (in the central slice) is shifted closer to zero for reconstructions with our filters, and the tail of the histogram is shorter. This highlights the fact that the \emph{maximum} standard deviation between reconstructions with our filters is smaller than the maximum standard deviation in reconstructions with the Shepp-Logan or Ram-Lak filters.

\subsubsection{Dependence of filters on noise and sparse angular sampling} 
\begin{figure}
    \includegraphics[scale=0.2]{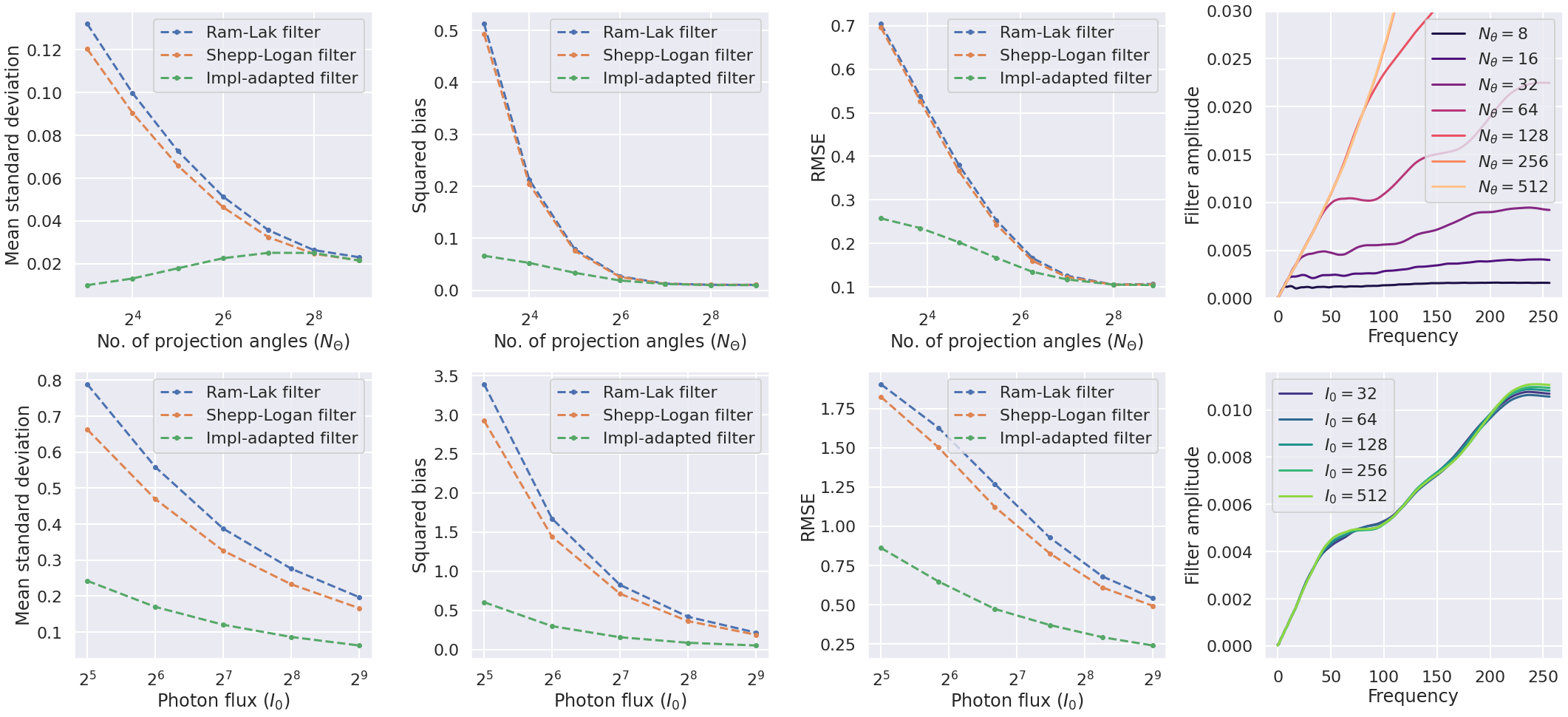}
    \caption{Implementation-adapted filters for noisy and sparsely sampled data. (\textit{top, left to right}) Mean standard deviations $\bar{\sigma}^S$ for slice $S = 128$ as a function of the number of projection angles $N_\theta$, mean value of the squared bias, mean value of RMSE with respect to the ground truth slice, and optimised filters in Fourier space. (\textit{bottom, left to right}) Mean standard deviations in $S=128$ as a function of photon flux $I_0$ (higher values of $I_0$ correspond to lower noise levels) using $N_\theta = 64$, mean value of the squared bias, mean value of RMSE with respect to the ground truth slice, and optimised filters in Fourier space.}
    \label{fig:foam_exp}
\end{figure}

We consider the effect of noise and sparse sampling on our filters. For the central slice of the foam phantom shown in Fig.~\ref{fig:foam}, we generated data by varying the number of projection angles $N_\theta$ and the photon flux $I_0$. For each of these settings, we computed the mean standard deviation \eqref{eq:mean_std_dev} between reconstruction slices. Our results are shown in Fig.~\ref{fig:foam_exp}. For all noise and angular sampling settings, the mean standard deviation in the slice was reduced by using implementation-adapted filters, with the difference being particularly prominent for noisy and smaller angular sampling settings. Shepp-Logan filter reconstructions had smaller mean standard deviation compared with Ram-Lak filter reconstructions, except in situations where many angles ($N_\theta\geq 256$) were used. In the high angle regime, reconstructions using the Ram-Lak filter have a relatively small number of artefacts and improvements due to filter optimisation are modest. 

We also quantified the mean squared bias and the mean RMSE with respect to the ground truth for this slice. From these plots, we observe that reconstructions using implementation-adapted filters have lower mean squared bias and mean RMSE compared with those for reconstructions with standard filters. High noise (low $I_0$) and sparse angular sampling settings result in an increase in bias and RMSE for all filter types. However, the increase is sharper for the Shepp-Logan and Ram-Lak filters than for our implementation-adapted filters. For every noise setting, the Ram-Lak filter results in the worst reconstructions in terms of bias and RMSE. Although both bias and RMSE increase as the number of projection angles is reduced in the noise-free setting, we observe a reduction in mean standard deviation for reconstructions using implementation-adapted filters. This suggests that in spite of a reduction in mean standard deviation due to effective suppression of high frequencies, the reconstructions produced by our implementation-adapted filters in this regime are incapable of mitigating the large number of low-angle artefacts. In effect, these settings show a limit where optimisation of a linear filter is not sufficient for good reconstructions, and intra-set homogeneity is achieved at the expense of an increase in bias and RMSE.

In addition, we also show the shapes of the filters (computed for the \texttt{strip} kernel in the ASTRA toolbox) as a function of noise and angular sampling. As the number of projection angles is increased, the shape of implementation-adapted filters approaches that of the ramp filter. In these regimes, reconstructions obtained using the Ram-Lak filter and the Shepp-Logan filter are nearly identical in terms of bias and RMSE. For different noise settings, the filters only vary at certain frequencies. It is possible that these frequencies are indicative of the main features in the foam phantom slice used.

\subsubsection{Variation of filters with projection data}
\begin{figure}
    \centering
    \includegraphics[scale=0.2]{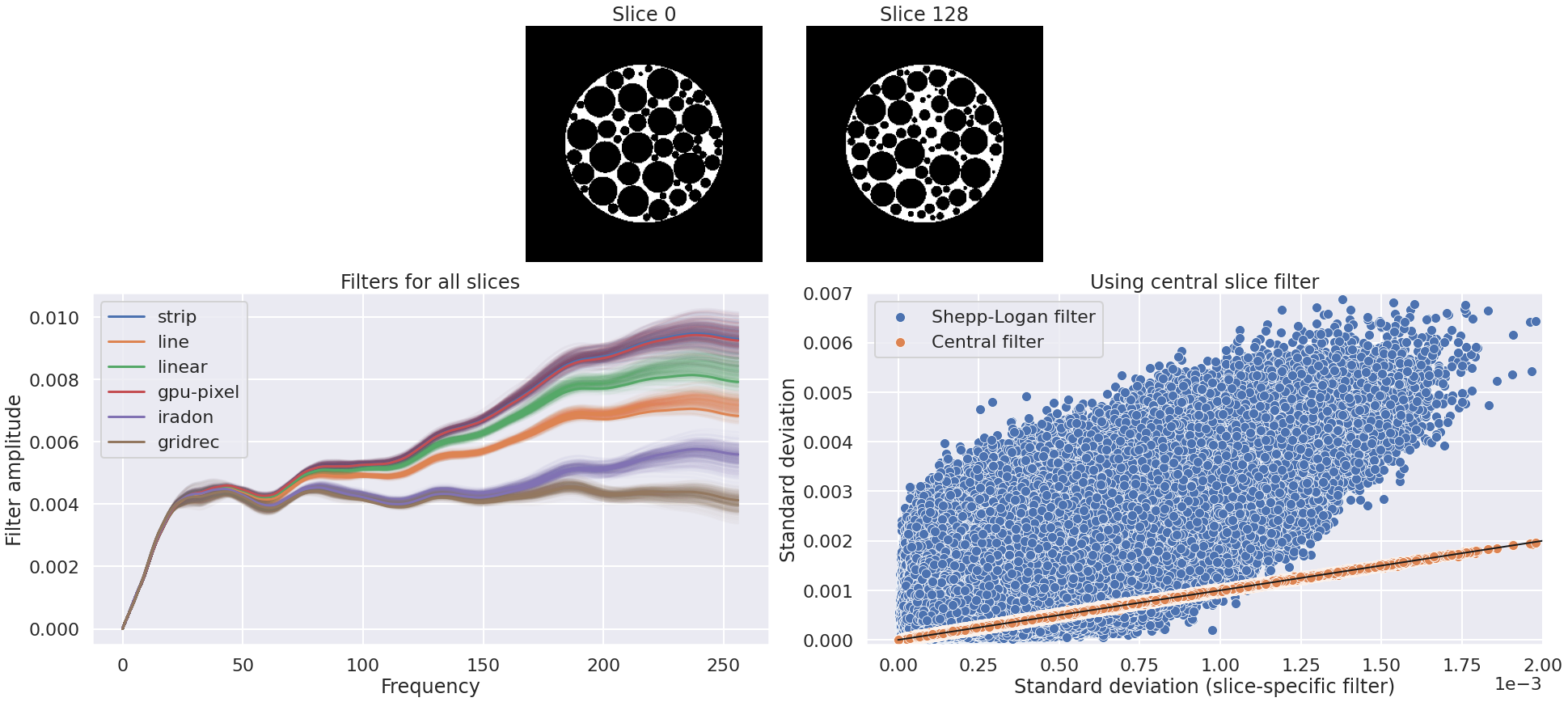}
    \caption{Variation of filters with projection data. (\textit{top}) Two slices of a simulated foam phantom with differences in features. (\textit{bottom left}) Implementation-adapted filters for all slices of the foam phantom (slice-specific filters). Central slice (slice no.~128) filters for each implementation are indicated with bold lines. (\textit{bottom right}) Scatter plot of pixelwise standard deviations $\sigma$ using slice-specific filters, the central slice filter and the Shepp-Logan filter. Standard deviations using the central slice filter are almost the same as those using slice-specific filters (orange dots). These points lie on a straight line (shown in black) with slope $\sim 1$ and intercept $ \sim 0$. In contrast, standard deviations using the Shepp-Logan filter are higher than those using slice-specific filters (blue dots) for most pixels.}
    \label{fig:filter_variability}
\end{figure}

In order to understand how our filters change with changes in the data, we computed filters for all slices of our simulated foam phantom. Two such slices are shown in Fig.~\ref{fig:filter_variability}. These slices, although visually similar, have different features. Implementation-adapted filters for all $256$ slices of the foam phantom are shown in Fig.~\ref{fig:filter_variability}. 

In order to study the applicability of the central slice filter to other slices, we performed the following experiment. First, we reconstructed all slices using the slice-specific filters, i.e.~filters that had been optimised for \emph{each individual slice} using different implementations. Next, we reconstructed all slices with the central slice filter. As a baseline, all slices were also reconstructed using the Shepp-Logan filter. Pixelwise standard deviations \eqref{eq:pixelwise_std} were computed for all pixels in the foam phantom volume for the three cases. The scatter plot in Fig.~\ref{fig:filter_variability} shows that the pixelwise standard deviations with the central slice filter are nearly the same as those with the slice-specific filters. In fact, these points lie on a line with slope nearly equal to one. This indicates that using the central slice filter results in an equivalent reduction in differences between reconstructions as slice-specific filters. In contrast, the pixelwise standard deviations using the Shepp-Logan filter are, for a majority of pixels, larger than those obtained using slice-specific filters. This suggests that, for a majority of pixels in the reconstruction volume, smaller values of standard deviation are observed after filter optimisation.

Our experiment suggests that using the central slice filter for all slices of the foam phantom results in an equivalent reduction in standard deviation as slice-specific filters. This paves the way to fast application of such filters in a real dataset. An implementation-adapted filter computed for one slice of such a dataset could be reused with all other slices with no additional computational cost, just like any of the standard filters in a software package.

\subsubsection{Reduction in differences after thresholding}
\begin{figure}
    \includegraphics[scale=0.3]{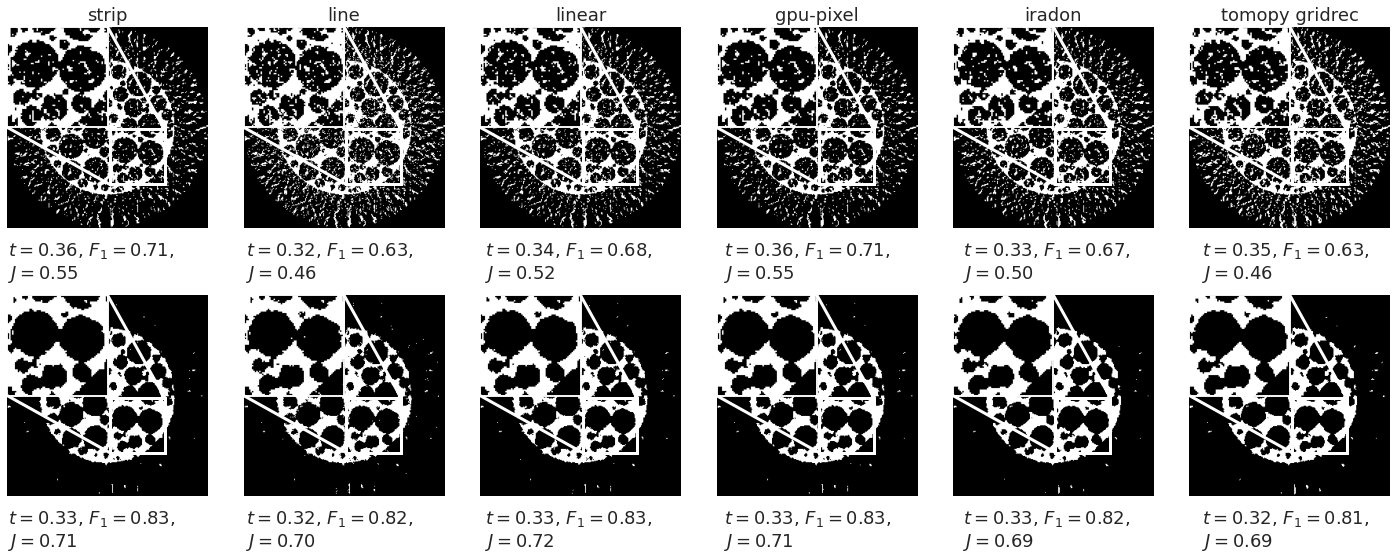}
    \caption{Differences after thresholding using Otsu's method. Reconstructions shown in Fig.~\ref{fig:foam} were used as input to the thresholding routine. (\textit{top row}) Thresholded reconstructions obtained using different backprojector implementations and the Shepp-Logan filter. Corresponding Otsu thresholds $t$, $F_1$ scores and Jaccard indices are given for each image. (\textit{bottom row}) Thresholded reconstructions obtained using implementation-adapted filters.}
    \label{fig:seg_conn}
\end{figure}

We investigated the effect of our filters on the results of a simple post-processing step. We reconstructed data ($N_\theta = 32$, no noise) from the central slice of the foam phantom and used Otsu's method in \texttt{scikit-image} \cite{van2014scikit} to threshold reconstruction slices from different implementations. In Fig.~\ref{fig:seg_conn}, we show two sets of thresholded reconstructions, one obtained using the Shepp-Logan filter and the other obtained using our implementation-adapted filters. We show values for the Otsu threshold $t$, the $F_1$ score with respect to the ground truth slice and the Jaccard index in the figure. We used routines in \texttt{scikit-learn} \cite{scikit-learn} to compute all segmentation metrics. For the set of Shepp-Logan filter reconstructions, the ranges of threshold values (0.32-0.36), $F_1$ scores (0.63-0.71) and Jaccard indices (0.46-0.55) were larger than the corresponding ranges for the implementation-adapted filter reconstructions. For the latter set, the Otsu threshold varied between 0.32 and 0.33 for all reconstructions. The $F_1$ scores were between 0.81 and 0.83, and the Jaccard indices were in the range of 0.69-0.72. Upon visual inspection of the zoomed-in insets we find greater differences between thresholded reconstructions in the set of Shepp-Logan filter reconstructions. These results suggest that post-processing steps such as segmentation may be rendered more reproducible and amenable to automation if reconstructions are obtained using implementation-adapted filters. 

\subsubsection{Optimising to a reference reconstruction}
Although we focus on filter optimisation in sinogram space in this paper, a related optimisation problem is one where reconstruction results from different implementations are optimised to a reference reconstruction. This type of optimisation might be useful when the result of one specific implementation is preferred due to its superior accuracy and when the exact settings used with this algorithm are unknown. 

In some cases, high-quality reconstructions might be computed with an unknown (possibly in-house) software package during the experiment by expert beamline scientists. When users reconstruct this data later at their home institutes, it might not be possible to use the same software packages with identical settings. Our approach would enable users to reduce the difference between their reconstructions and the high-quality reference reconstructions.

Optimisation in reconstruction space can be performed by modifying the objective in \eqref{eq:filter_opt}:
\begin{equation}
    \bm{h}_I^\ast = \arg \min_{\bm{h}} \|\bm{r}_{\text{ref}} -  \bm{r}_I(\bm{h}, \bm{p})\|^2_2,
\end{equation}
where $\bm{r}_{\text{ref}}$ is the reference reconstruction. 

\begin{figure}
    \centering
    \includegraphics[scale=0.2]{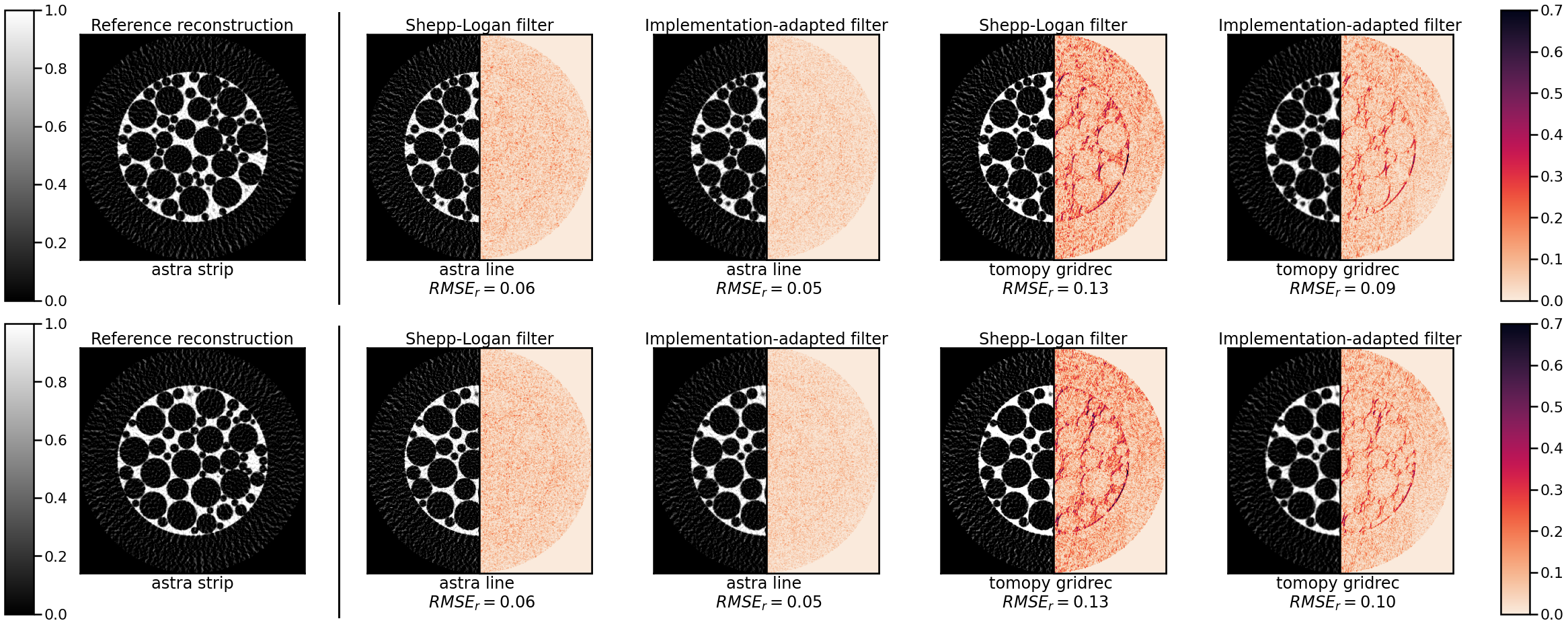}
    \caption{Filter optimisation using a reference reconstruction. (\textit{top row}) Filters optimised to a \texttt{strip} kernel reconstruction (\textit{top row, left}). (\textit{top row}) Reconstructions before and after filter optimisation using the ASTRA \texttt{line} kernel and Gridrec. Right half of each image shows absolute difference with the reference reconstruction. RMSE values with respect to the reference are also shown.
    (\textit{bottom row}) Reconstructions of a different (test) slice using the filters obtained for the slice in the top row. Pixelwise absolute difference and RMSE using implementation-adapted filters are smaller in both cases.}
    \label{fig:opt_to_ref_reco}
\end{figure}

To illustrate filter optimisation in reconstruction space, we performed the following experiment. Using the \texttt{strip} kernel reconstruction (with the Shepp-Logan filter) as a reference, we computed optimised filters for two other implementations (ASTRA \texttt{line} kernel and TomoPy Gridrec) for reconstructing the central slice of the foam phantom. Subsequently, we reconstructed the sinogram with the Shepp-Logan filter and our filters. These reconstructions are shown in the top row of Fig.~\ref{fig:opt_to_ref_reco}. To quantify similarity with the reference reconstruction, we computed the pixelwise absolute difference between each reconstruction and the reference as well as the RMSE using the reference as ground truth, which we denote as $\text{RMSE}_r$. For both \texttt{line} and Gridrec backprojectors, optimising the filter to a reference reconstruction reduced the $\text{RMSE}_r$ and absolute difference. 
As a further test, we applied the filters computed for this slice to a different slice of the foam phantom, which did not have any overlaps with the slice used to compute the filters. For this test slice, we again observed the reduction in $\text{RMSE}_r$ and absolute error, suggesting that our filters were able to bring the resulting reconstructions closer to the reference reconstruction.

\subsection{Round-Robin data}
\begin{figure}
\includegraphics[scale=0.19]{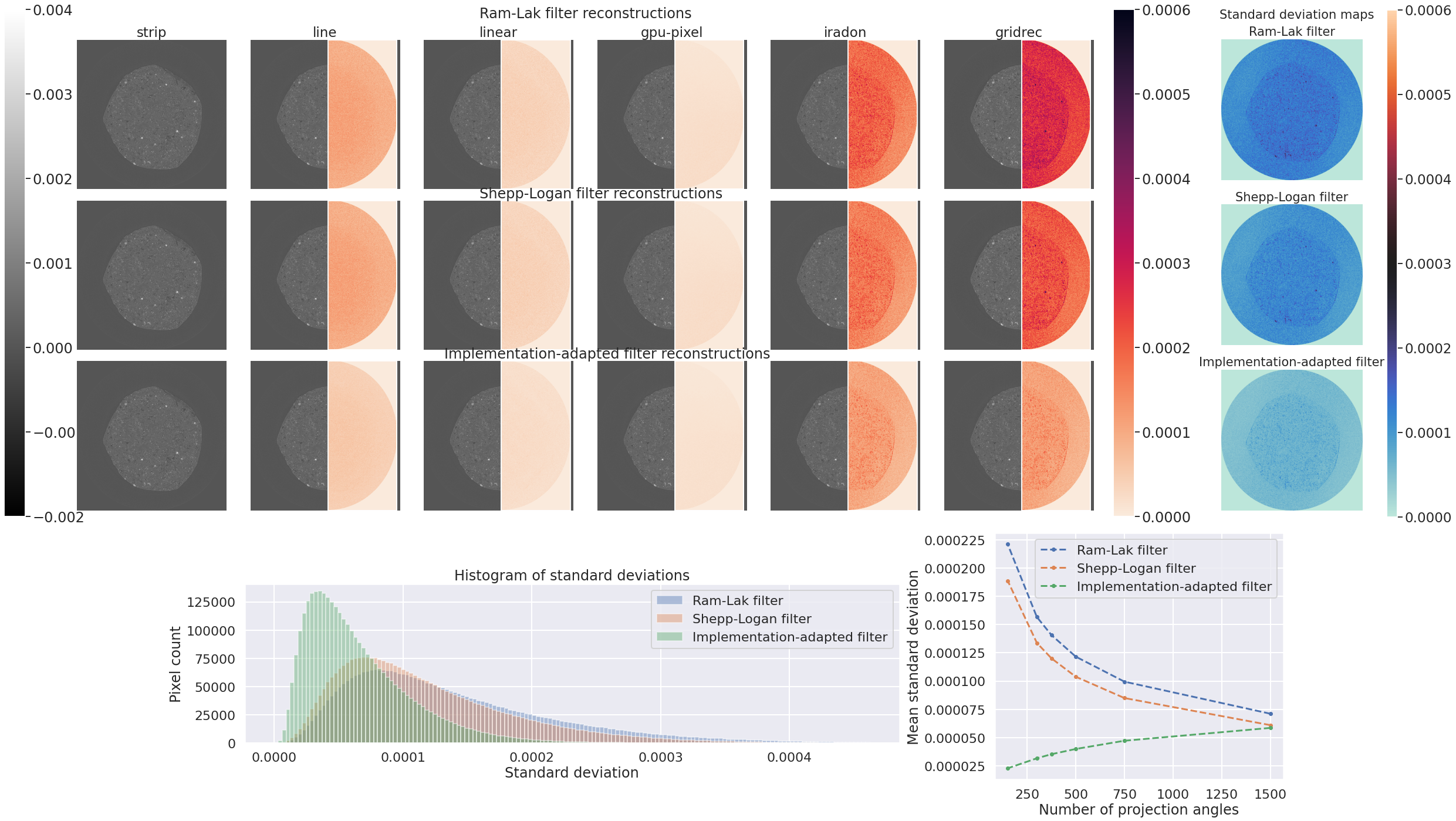}
    \caption{Reduction in differences between reconstructions of the Round-Robin dataset (slice no.~896). (\textit{top three rows}) Slice reconstructions using different implementations. Reconstructions were performed by discarding every second projection from the full dataset. The right half of the images show absolute differences with the corresponding \texttt{strip} kernel reconstruction in each set. The rightmost column shows pixelwise standard deviations in each set.
    (\textit{bottom row, left})
    Histograms of standard deviation for all three types of filters. (\textit{right}) Mean standard deviations $\bar{\sigma}^S$ in slice $S = 896$ for different numbers of projection angles.}
    \label{fig:roundrobin}
\end{figure}

Fig.~\ref{fig:roundrobin} shows the results of our method on the central slice (slice no.~896) of the Round-Robin dataset N1. These reconstructions were performed by discarding every second projection from the entire dataset. From the heatmaps of absolute difference with respect to the \texttt{strip} kernel reconstruction, we observe that intra-set differences are reduced by using implementation-adapted filters. This is further shown by the pixelwise standard deviation maps. Standard deviations between reconstructions using the Ram-Lak and Shepp-Logan filters are larger than those between reconstructions using implementation-adapted filters. Similar to the distributions in Fig.~\ref{fig:foam}, we see that our implementation-adapted filters are able to shift the mode of the histogram of standard deviations towards zero and to reduce the number of large standard deviations in the slice. We also observe that the Ram-Lak filter reconstructions show higher standard deviations than the Shepp-Logan filter reconstructions.

We also studied the effect of the number of projections used on the mean standard deviation \eqref{eq:mean_std_dev} in this slice. To do this, we performed experiments with the whole dataset and also with parts of the data, where every $2, 3, 4, 5 \text { and } 10$ projections were discarded. For each instance, the data were reconstructed using the Ram-Lak filter, the Shepp-Logan filter and our implementation-adapted filters. The plot of mean standard deviations is shown in Fig.~\ref{fig:roundrobin}. For all projection numbers, filter optimisation reduced the mean standard deviation in the slice. The difference was smaller for higher projection numbers, indicating that our filters are especially useful in improving reproducibility of reconstructions when the number of projection angles is small. In practice, data along few angles may be acquired to reduce the X-ray dose on a sample or to speed up acquisition when the sample is evolving over time.

\section{Discussion}\label{discussion}

In this paper, we presented a method to improve the reproducibility of reconstructions in the synchrotron pipeline. Our method uses an optimisation problem over filters to reduce differences between reconstructions from various software implementations of commonly-used algorithms.

The objective function that was used in our optimisation problem was the $\ell^2$-distance between the forward projection of the obtained reconstruction and the given projection data. This choice was motivated by the fact that ground truth reconstructions are generally not available in real-world experiments. However, it is possible to formulate a similar (and related) problem in reconstruction space, by using the $\ell^2$-distance between the reconstruction from a given software package and a reference reconstruction as the objective to be minimised. The solution to such an optimisation procedure is a shift-invariant blurring kernel in reconstruction space. The implementation-adapted filters presented in this paper can thus be viewed as a linear transformation of the projection data that results in an automatic selection of shift-invariant blurring of reconstructions.

Our work here can be extended to optimise other pre-processing and post-processing steps in the synchrotron pipeline. An important example is phase retrieval, which can be formulated in terms of a filtering operation \cite{paganin2002simultaneous}. This filter can be optimised similarly in order to improve reproducibility.

\begin{figure}
    \centering
    \includegraphics[scale=0.4]{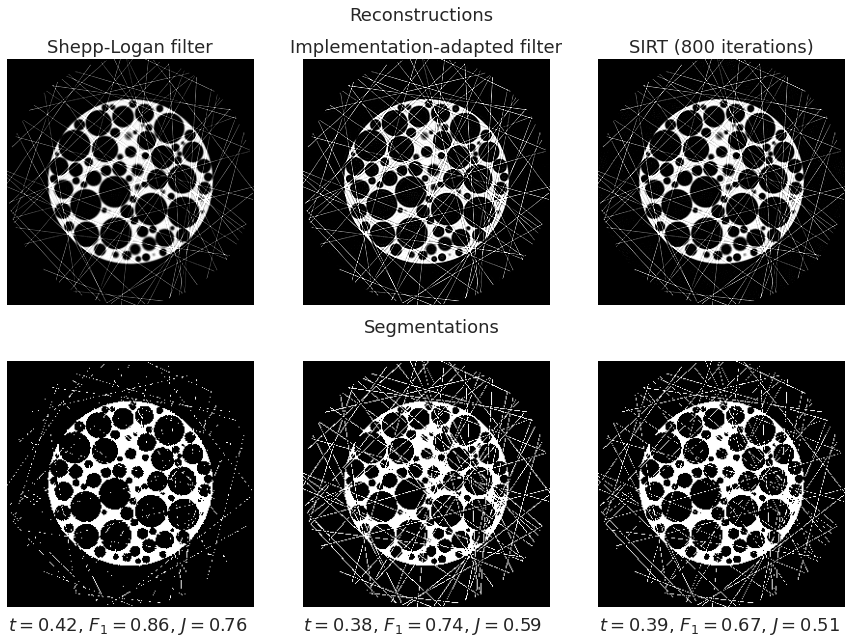}
    \caption{Reconstructions of data corrupted with zingers showing an example where the Shepp-Logan filter reconstruction and corresponding segmentation are better than those using an implementation-adapted filter or an iterative method (SIRT). (\textit{top row}) Reconstructions of data from slice 128 ($N_\theta=512$, no noise) corrupted with zingers. Zingers are more prominent in the reconstruction using an implementation-adapted filter and in the SIRT reconstruction (after 800 iterations). (\textit{bottom row}) Segmentations using Otsu's method of all three reconstructions. The Otsu threshold, $F_1$ score and Jaccard index for each image is given below.}
    \label{fig:zingers}
\end{figure}

One limitation of our method is that we optimise to the data available. This optimisation can lead to undesired solutions in the presence of outliers in the data, such as zingers or ring artefacts. Reconstructions of data corrupted with zingers (randomly placed very bright pixels in the sinogram) are shown in Fig.~\ref{fig:zingers}. In this example we see that the FBP reconstruction using the ASTRA \texttt{strip} kernel and the Shepp-Logan filter shows less prominent zingers than the reconstruction using an implementation-adapted filter. This is because the optimised filter preserves the zingers in the data whereas the unoptimised FBP reconstruction is independent of them. Other methods, such as the simultaneous iterative reconstruction technique (SIRT), which iteratively minimise the data misfit also give similar, poor reconstructions. One way to improve iterative reconstruction methods is to use regularisation, which can be achieved either by early stopping or by the inclusion of an explicit regularisation term in the objective function to be minimised. Analogous techniques can be used for our filter optimisation problem \eqref{eq:filter_opt} to ensure greater robustness to outliers. 

Although we have demonstrated the reusability of our filters for similar data, these filters are dependent on the noise statistics and angular sampling in the acquired projections. One way to improve the generalisability of filters would be to simultaneously optimise to more than one dataset. This idea has been explored in \cite{pelt2013fast, lagerwerf2020computationally} using shallow neural networks.

Another promising direction is provided by deep learning-based methods, which have been applied to improve tomographic image reconstruction in a number of ways \cite{arridge2019solving}. Supervised deep learning approaches can be used to learn a (non-linear) mapping from input reconstructions to a reference reconstruction. However, such approaches generally require large amounts of paired training data (input and reference reconstructions). When insufficient training pairs are available, various unsupervised approaches, such as the Deep Image Prior method proposed in \cite{ulyanov2018deep}, are more suitable. For a quantitative comparison of various popular deep learning-based reconstruction methods, we refer the reader to \cite{leuschner2021quantitative}.

Apart from software solutions for image reconstruction, which have been the focus of this paper, improving reproducibility throughout the synchrotron pipeline requires hardware adjustments to the blocks in Fig~\ref{fig:pipeline}. Scanning the same sample twice under the same experimental conditions leads to small fluctuations in the data due to stochastic noise and drifts during the scanning process. In addition, beam-sensitive samples might deform due to irradiation. Such changes lead to differences in reconstructions that are similar to the differences due to software implementations, albeit less structured than those shown in Fig.~\ref{fig:diffs}. To improve hardware reproducibility, controlled phantom experiments might be performed to address differences in data acquisition. Finally, software and hardware solutions can be effectively linked by using approaches like reinforcement learning for experimental design and control \cite{recht2019tour, kain2020sample}. Such creative solutions might provide an efficient way for synchrotron users to perform reproducible experiments in the future.

\section{Conclusion}\label{conclusion}
In this paper, we proposed a filter optimisation method to improve reproducibility of tomographic reconstructions at synchrotrons. These implementation-adapted filters can be computed for any black-box software implementation by using only evaluations of the corresponding reconstruction routine. We numerically demonstrated the properties of and use cases for such filters. In both real and simulated data, our implementation-adapted filters reduced the standard deviation between reconstructions from various software implementations of reconstruction algorithms. The reduction in standard deviation was especially evident when the data were noisy or sparsely sampled. 

Our filter optimisation technique can be used to reduce the effect of differences in discretisation and interpolation in commonly-used software packages and is a key building block towards improving reproducibility throughout the synchrotron pipeline. We make available the open-source Python code for our method, allowing synchrotron users to obtain reconstructions that are more comparable and reproducible.

\ack{\textbf{Funding Information} }
P.S.G.~would like to acknowledge the financial support of the Marie Skłodowska-Curie Innovative Training Network MUMMERING (grant agreement no. 765604).  D.M.P.~is financially supported by The Netherlands Organization for Scientific Research (NWO), project number 016.Veni.192.235. F.d.C and D.G.'s work was supported by the U.S. Department of Energy, Office of Science and Technology, under contract DE-AC02-06CH11357.


\referencelist{iucr}






@inproceedings{palenstijn2013astra,
  title={The ASTRA tomography toolbox},
  author={Palenstijn, Willem Jan and Batenburg, K Joost and Sijbers, Jan},
  booktitle={13th International Conference on Computational and Mathematical Methods in Science and Engineering, CMMSE},
  volume={2013},
  pages={1139--1145},
  year={2013}
}

@article{van2014scikit,
  title={scikit-image: image processing in Python},
  author={Van der Walt, Stefan and Sch{\"o}nberger, Johannes L and Nunez-Iglesias, Juan and Boulogne, Fran{\c{c}}ois and Warner, Joshua D and Yager, Neil and Gouillart, Emmanuelle and Yu, Tony},
  journal={PeerJ},
  volume={2},
  pages={e453},
  year={2014},
  publisher={PeerJ Inc.}
}

@article{gursoy2014tomopy,
  title={TomoPy: a framework for the analysis of synchrotron tomographic data},
  author={G{\"u}rsoy, Doga and De Carlo, Francesco and Xiao, Xianghui and Jacobsen, Chris},
  journal={Journal of synchrotron radiation},
  volume={21},
  number={5},
  pages={1188--1193},
  year={2014},
  publisher={International Union of Crystallography}
}

@article{de2018tomobank,
  title={TomoBank: a tomographic data repository for computational x-ray science},
  author={De Carlo, Francesco and G{\"u}rsoy, Do{\u{g}}a and Ching, Daniel J and Batenburg, K Joost and Ludwig, Wolfgang and Mancini, Lucia and Marone, Federica and Mokso, Rajmund and Pelt, Dani{\"e}l M and Sijbers, Jan and others},
  journal={Measurement Science and Technology},
  volume={29},
  number={3},
  pages={034004},
  year={2018},
  publisher={IOP Publishing}
}

@inproceedings{xu2006comparative,
  title={A comparative study of popular interpolation and integration methods for use in computed tomography},
  author={Xu, Fang and Mueller, Klaus},
  booktitle={3rd IEEE International Symposium on Biomedical Imaging: Nano to Macro, 2006.},
  pages={1252--1255},
  year={2006},
  organization={IEEE}
}

@article{marone2012regridding,
  title={Regridding reconstruction algorithm for real-time tomographic imaging},
  author={Marone, F and Stampanoni, M},
  journal={Journal of synchrotron radiation},
  volume={19},
  number={6},
  pages={1029--1037},
  year={2012},
  publisher={International Union of Crystallography}
}

@article{arcadu2016fast,
  title={Fast gridding projectors for analytical and iterative tomographic reconstruction of differential phase contrast data},
  author={Arcadu, Filippo and Stampanoni, Marco and Marone, Federica},
  journal={Optics Express},
  volume={24},
  number={13},
  pages={14748--14764},
  year={2016},
  publisher={Optical Society of America}
}

@article{pelt2014improving,
  title={Improving filtered backprojection reconstruction by data-dependent filtering},
  author={Pelt, Dani{\"e}l M and Batenburg, Kees Joost},
  journal={IEEE Transactions on Image Processing},
  volume={23},
  number={11},
  pages={4750--4762},
  year={2014},
  publisher={IEEE}
}

@inproceedings{dowd1999developments,
  title={Developments in synchrotron x-ray computed microtomography at the National Synchrotron Light Source},
  author={Dowd, Betsy A and Campbell, Graham H and Marr, Robert B and Nagarkar, Vivek V and Tipnis, Sameer V and Axe, Lisa and Siddons, D Peter},
  booktitle={Developments in X-ray Tomography II},
  volume={3772},
  pages={224--236},
  year={1999},
  organization={International Society for Optics and Photonics}
}

@article{pelt2018improving,
  title={Improving tomographic reconstruction from limited data using mixed-scale dense convolutional neural networks},
  author={Pelt, Dani{\"e}l M and Batenburg, Kees Joost and Sethian, James A},
  journal={Journal of Imaging},
  volume={4},
  number={11},
  pages={128},
  year={2018},
  publisher={Multidisciplinary Digital Publishing Institute}
}

@article{midgley2009electron,
  title={Electron tomography and holography in materials science},
  author={Midgley, Paul A and Dunin-Borkowski, Rafal E},
  journal={Nature materials},
  volume={8},
  number={4},
  pages={271--280},
  year={2009},
  publisher={Nature Publishing Group}
}

@article{rubin2014computed,
  title={Computed tomography: revolutionizing the practice of medicine for 40 years},
  author={Rubin, Geoffrey D},
  journal={Radiology},
  volume={273},
  number={2S},
  pages={S45--S74},
  year={2014},
  publisher={Radiological Society of North America}
}

@article{fusseis2014brief,
  title={A brief guide to synchrotron radiation-based microtomography in (structural) geology and rock mechanics},
  author={Fusseis, Florian and Xiao, Xianghui and Schrank, Christoph and De Carlo, F},
  journal={Journal of Structural Geology},
  volume={65},
  pages={1--16},
  year={2014},
  publisher={Elsevier}
}

@article{luo2018cracking,
  title={Cracking evolution behaviors of lightweight materials based on in situ synchrotron X-ray tomography: A review},
  author={Luo, Y and Wu, SC and Hu, YN and Fu, YN},
  journal={Frontiers of Mechanical Engineering},
  volume={13},
  number={4},
  pages={461--481},
  year={2018},
  publisher={Springer}
}

@misc{kak2002principles,
  title={Principles of computerized tomographic imaging},
  author={Kak, Avinash C and Slaney, Malcolm and Wang, Ge},
  year={2002},
  publisher={Wiley Online Library}
}

@incollection{buzug2011computed,
  title={Computed tomography},
  author={Buzug, Thorsten M},
  booktitle={Springer Handbook of Medical Technology},
  pages={311--342},
  year={2011},
  publisher={Springer}
}

@book{natterer2001mathematics,
  title={The mathematics of computerized tomography},
  author={Natterer, Frank},
  year={2001},
  publisher={SIAM}
}

@article{thompson1984computed,
  title={Computed tomography using synchrotron radiation},
  author={Thompson, AC and Llacer, J and Finman, L Campbell and Hughes, EB and Otis, JN and Wilson, S and Zeman, HD},
  journal={Nuclear Instruments and Methods in Physics Research},
  volume={222},
  number={1-2},
  pages={319--323},
  year={1984},
  publisher={Elsevier}
}

@inproceedings{de2006x,
  title={X-ray tomography system, automation, and remote access at beamline 2-BM of the Advanced Photon Source},
  author={De Carlo, Francesco and Xiao, Xianghui and Tieman, Brian},
  booktitle={Developments in X-ray Tomography V},
  volume={6318},
  pages={63180K},
  year={2006},
  organization={International Society for Optics and Photonics}
}

@book{stock2019microcomputed,
  title={Microcomputed tomography: methodology and applications},
  author={Stock, Stuart R},
  year={2019},
  publisher={CRC press}
}

@article{hintermuller2010image,
  title={Image processing pipeline for synchrotron-radiation-based tomographic microscopy},
  author={Hinterm{\"u}ller, C and Marone, F and Isenegger, A and Stampanoni, M},
  journal={Journal of synchrotron radiation},
  volume={17},
  number={4},
  pages={550--559},
  year={2010},
  publisher={International Union of Crystallography}
}

@article{paganin2002simultaneous,
  title={Simultaneous phase and amplitude extraction from a single defocused image of a homogeneous object},
  author={Paganin, David and Mayo, Sheridan C and Gureyev, Tim E and Miller, Peter R and Wilkins, Steve W},
  journal={Journal of microscopy},
  volume={206},
  number={1},
  pages={33--40},
  year={2002},
  publisher={Wiley Online Library}
}

@article{yang2017convolutional,
  title={A convolutional neural network approach to calibrating the rotation axis for X-ray computed tomography},
  author={Yang, Xiaogang and De Carlo, Francesco and Phatak, Charudatta and G{\"u}rsoy, Dogˇa},
  journal={Journal of Synchrotron Radiation},
  volume={24},
  number={2},
  pages={469--475},
  year={2017},
  publisher={International Union of Crystallography}
}

@article{massimi2018improved,
  title={An improved ring removal procedure for in-line x-ray phase contrast tomography},
  author={Massimi, Lorenzo and Brun, Francesco and Fratini, Michela and Bukreeva, Inna and Cedola, Alessia},
  journal={Physics in Medicine \& Biology},
  volume={63},
  number={4},
  pages={045007},
  year={2018},
  publisher={IOP Publishing}
}

@article{pelt2016integration,
  title={Integration of TomoPy and the ASTRA toolbox for advanced processing and reconstruction of tomographic synchrotron data},
  author={Pelt, Dani{\"e}l M and G{\"u}rsoy, Doga and Palenstijn, Willem Jan and Sijbers, Jan and De Carlo, Francesco and Batenburg, Kees Joost},
  journal={Journal of synchrotron radiation},
  volume={23},
  number={3},
  pages={842--849},
  year={2016},
  publisher={International Union of Crystallography}
}

@article{buhrer2020unveiling,
  title={Unveiling water dynamics in fuel cells from time-resolved tomographic microscopy data},
  author={B{\"u}hrer, Minna and Xu, Hong and Eller, Jens and Sijbers, Jan and Stampanoni, Marco and Marone, Federica},
  journal={Scientific Reports},
  volume={10},
  number={1},
  pages={1--15},
  year={2020},
  publisher={Nature Publishing Group}
}

@article{salome1999synchrotron,
  title={A synchrotron radiation microtomography system for the analysis of trabecular bone samples},
  author={Salom{\'e}, Murielle and Peyrin, Fran{\c{c}}oise and Cloetens, Peter and Odet, Christophe and Laval-Jeantet, Anne-Marie and Baruchel, Jos{\'e} and Spanne, Per},
  journal={Medical Physics},
  volume={26},
  number={10},
  pages={2194--2204},
  year={1999},
  publisher={Wiley Online Library}
}

@article{lagerwerf2020computationally,
  title={A computationally efficient reconstruction algorithm for circular cone-beam computed tomography using shallow neural networks},
  author={Lagerwerf, Marinus J and Pelt, Dani{\"e}l M and Palenstijn, Willem Jan and Batenburg, Kees Joost},
  journal={Journal of Imaging},
  volume={6},
  number={12},
  pages={135},
  year={2020},
  publisher={Multidisciplinary Digital Publishing Institute}
}

@article{pelt2013fast,
  title={Fast tomographic reconstruction from limited data using artificial neural networks},
  author={Pelt, Daniel Maria and Batenburg, Kees Joost},
  journal={IEEE Transactions on Image Processing},
  volume={22},
  number={12},
  pages={5238--5251},
  year={2013},
  publisher={IEEE}
}

@inproceedings{ulyanov2018deep,
  title={Deep image prior},
  author={Ulyanov, Dmitry and Vedaldi, Andrea and Lempitsky, Victor},
  booktitle={Proceedings of the IEEE conference on computer vision and pattern recognition},
  pages={9446--9454},
  year={2018}
}

@article{kanitpanyacharoen2013comparative,
  title={A comparative study of X-ray tomographic microscopy on shales at different synchrotron facilities: ALS, APS and SLS},
  author={Kanitpanyacharoen, Waruntorn and Parkinson, Dilworth Y and De Carlo, Francesco and Marone, Federica and Stampanoni, Marco and Mokso, Rajmund and MacDowell, Alastair and Wenk, H-R},
  journal={Journal of synchrotron radiation},
  volume={20},
  number={1},
  pages={172--180},
  year={2013},
  publisher={International Union of Crystallography}
}

@incollection{pchansenbook2021,
  author      = {Batenburg, Kees Joost and Hansen, Per Christian and Jorgensen, Jakob Sauer},
  title       = {Discretization Models and the System Matrix},
  editor      = {Hansen, Per Christian and Jorgensen, Jakob Sauer and Lionheart, William R. B.},
  booktitle   = {Scientific Computing for Computed Tomography},
  publisher   = {in press},
  year        = 2021,
  chapter     = 8,
}

@Article{         harris2020array,
 title         = {Array programming with {NumPy}},
 author        = {Charles R. Harris and K. Jarrod Millman and St{'{e}}fan J.
                 van der Walt and Ralf Gommers and Pauli Virtanen and David
                 Cournapeau and Eric Wieser and Julian Taylor and Sebastian
                 Berg and Nathaniel J. Smith and Robert Kern and Matti Picus
                 and Stephan Hoyer and Marten H. van Kerkwijk and Matthew
                 Brett and Allan Haldane and Jaime Fern{'{a}}ndez del
                 R{'{\i}}o and Mark Wiebe and Pearu Peterson and Pierre
                 G{'{e}}rard-Marchant and Kevin Sheppard and Tyler Reddy and
                 Warren Weckesser and Hameer Abbasi and Christoph Gohlke and
                 Travis E. Oliphant},
 year          = {2020},
 month         = sep,
 journal       = {Nature},
 volume        = {585},
 number        = {7825},
 pages         = {357--362},
 doi           = {10.1038/s41586-020-2649-2},
 publisher     = {Springer Science and Business Media {LLC}},
 url           = {https://doi.org/10.1038/s41586-020-2649-2}
}

@article{arridge2019solving,
  title={Solving inverse problems using data-driven models},
  author={Arridge, Simon and Maass, Peter and {\"O}ktem, Ozan and Sch{\"o}nlieb, Carola-Bibiane},
  journal={Acta Numerica},
  volume={28},
  pages={1--174},
  year={2019},
  publisher={Cambridge University Press}
}

@Article{leuschner2021quantitative,
AUTHOR = {Leuschner, Johannes and Schmidt, Maximilian and Ganguly, Poulami Somanya and Andriiashen, Vladyslav and Coban, Sophia Bethany and Denker, Alexander and Bauer, Dominik and Hadjifaradji, Amir and Batenburg, Kees Joost and Maass, Peter and Eijnatten, Maureen van},
TITLE = {Quantitative Comparison of Deep Learning-Based Image Reconstruction Methods for Low-Dose and Sparse-Angle CT Applications},
JOURNAL = {Journal of Imaging},
VOLUME = {7},
YEAR = {2021},
NUMBER = {3},
ARTICLE-NUMBER = {44},
URL = {https://www.mdpi.com/2313-433X/7/3/44},
ISSN = {2313-433X},
DOI = {10.3390/jimaging7030044}
}

@article{kain2020sample,
  title={Sample-efficient reinforcement learning for CERN accelerator control},
  author={Kain, Verena and Hirlander, Simon and Goddard, Brennan and Velotti, Francesco Maria and Della Porta, Giovanni Zevi and Bruchon, Niky and Valentino, Gianluca},
  journal={Physical Review Accelerators and Beams},
  volume={23},
  number={12},
  pages={124801},
  year={2020},
  publisher={APS}
}

@article{recht2019tour,
  title={A tour of reinforcement learning: The view from continuous control},
  author={Recht, Benjamin},
  journal={Annual Review of Control, Robotics, and Autonomous Systems},
  volume={2},
  pages={253--279},
  year={2019},
  publisher={Annual Reviews}
}

@article{lagerwerf2020automated,
  title={Automated FDK-filter selection for Cone-beam CT in research environments},
  author={Lagerwerf, Marinus J and Palenstijn, Willem Jan and Kohr, Holger and Batenburg, Kees Joost},
  journal={IEEE Transactions on Computational Imaging},
  volume={6},
  pages={739--748},
  year={2020},
  publisher={IEEE}
}

@article{otsu1979threshold,
  title={A threshold selection method from gray-level histograms},
  author={Otsu, Nobuyuki},
  journal={IEEE transactions on systems, man, and cybernetics},
  volume={9},
  number={1},
  pages={62--66},
  year={1979},
  publisher={IEEE}
}

@article{scikit-learn,
 title={Scikit-learn: Machine Learning in {P}ython},
 author={Pedregosa, F. and Varoquaux, G. and Gramfort, A. and Michel, V.
         and Thirion, B. and Grisel, O. and Blondel, M. and Prettenhofer, P.
         and Weiss, R. and Dubourg, V. and Vanderplas, J. and Passos, A. and
         Cournapeau, D. and Brucher, M. and Perrot, M. and Duchesnay, E.},
 journal={Journal of Machine Learning Research},
 volume={12},
 pages={2825--2830},
 year={2011}
}
\end{document}